\documentclass{article}
\usepackage[utf8]{inputenc}

\usepackage{xcolor}
\usepackage[centering,marginparwidth=0.5cm]{geometry}  
\usepackage{setspace}
\usepackage[title]{appendix}
\usepackage[round]{natbib}
\usepackage{subcaption}
\usepackage{amssymb}
\usepackage{amsmath}
\usepackage{bbm}
\usepackage{mathtools}
\usepackage{multicol}
\usepackage{graphicx}
\usepackage{array}
\usepackage{bm}
\usepackage{soul}
\usepackage{dirtytalk}
\usepackage{makecell}
\usepackage{caption}
\usepackage[noend]{algpseudocode}



\usepackage{algorithm,algpseudocode,caption}

\makeatletter
\newcommand\fs@boxedtopcap{\def\@fs@cfont{\bfseries}\let\@fs@capt\floatc@plain
	\def\@fs@pre{\setbox\@currbox\vbox{\hbadness10000
			\moveleft3.4pt\vbox{\advance\hsize by6.8pt
				\hrule \hbox to\hsize{\vrule\kern3pt
					\vbox{\kern3pt\box\@currbox\kern3pt}\kern3pt\vrule}\hrule}}}%
	\def\@fs@mid{\kern2pt}%
	\def\@fs@post{}\let\@fs@iftopcapt\iftrue}
\makeatother

\floatstyle{boxedtopcap}
\restylefloat{algorithm}

\usepackage{threeparttable}

\usepackage{enumitem}
\setlist[enumerate]{noitemsep, topsep=2pt}
\setlist[itemize]{noitemsep, topsep=2pt}

\usepackage{extarrows}

\usepackage{anysize}
\usepackage{latexsym}
\usepackage{amsthm}
\usepackage{epsfig}

\usepackage[T1]{fontenc}
\usepackage{mathrsfs}
\usepackage{float}
\usepackage{setspace}
\usepackage{color}
\definecolor{lawngreen}{RGB}{0,250,154}

\usepackage{titletoc}

\usepackage{color}
\usepackage{dsfont}
\usepackage{booktabs}
\definecolor{darkblue}{rgb}{0.0,0.0,0.5}
\definecolor{winered}{rgb}{0.5,0,0}
\definecolor{deeppink}{RGB}{255,20,147}
\usepackage[bookmarks,colorlinks,breaklinks]{hyperref}  
\hypersetup{linkcolor=winered,citecolor=winered,filecolor=winered,urlcolor=winered} 
 
\allowdisplaybreaks[1]

\graphicspath{{./Figures/}}

\usepackage{etoolbox}
\usepackage{lipsum}
\makeatletter
\patchcmd{\@addmarginpar}{\ifodd\c@page}{\ifodd\c@page\@tempcnta\m@ne}{}{}
\makeatother

\usepackage{thmtools}

\theoremstyle{definition}
\newtheorem{theorem}{Theorem}

\def\phcomments{1}
\newcounter{note}[section]
\renewcommand{\thenote}{\thesection.\arabic{note}}
\newcommand{\noteR}[2]{\refstepcounter{note}\marginpar{\tiny\bf \textcolor{red}{#1~\thenote}}
$\ll${\sf \textcolor{red}{#1~\thenote:}} 
{\color{red}{#2}}$\gg$}

\newcommand{\noteP}[2]{\refstepcounter{note}\marginpar{\tiny\bf \textcolor{deeppink}{#1~\thenote}}
$\ll${\sf \textcolor{deeppink}{#1~\thenote:}} 
{\color{deeppink}{#2}}$\gg$}

\newcommand{\noteO}[2]{\refstepcounter{note}\marginpar{\tiny\bf \textcolor{orange}{#1~\thenote}}
$\ll${\sf \textcolor{orange}{#1~\thenote:}} 
{\color{orange}{#2}}$\gg$}

\newcommand{\noteS}[2]{\refstepcounter{note}\marginpar{\tiny\bf \textcolor{blue}{#1~\thenote}}
$\ll${\sf \textcolor{orange}{#1~\thenote:}} 
{\color{blue}{#2}}$\gg$}

\ifnum\phcomments=1
\newcommand{\chris}[1]{\noteP{Chris}{#1}}
\newcommand{\zhenyu}[1]{\noteR{Zhenyu}{#1}}
\newcommand{\joline}[1]{\noteO{JU}{#1}}
\newcommand{\hs}[1]{\noteS{HS}{#1}}
\else
\newcommand{\chris}[1]{}
\newcommand{\zhenyu}[1]{}
\newcommand{\joline}[1]{}
\newcommand{\hs}[1]{}
\fi

\usepackage[nameinlink]{cleveref}
\crefname{assumption}{Assumption}{Assumptions}
\crefname{lemma}{Lemma}{Lemmas}
\crefname{theorem}{Theorem}{Theorems}
\crefname{corollary}{Corollary}{Corollaries}
\crefname{proposition}{Proposition}{Propositions}
\crefname{claim}{Claim}{Claims}
\crefname{procedure}{Procedure}{Procedures}
\crefname{algorithm}{Algorithm}{Algorithms}
\crefname{figure}{Figure}{Figures}
\crefname{remark}{Remark}{Remarks}
\crefname{section}{Section}{Sections}
\crefname{procedure}{Procedure}{Procedures}
\crefname{table}{Table}{Tables}
\crefname{equation}{}{}
\crefname{enumi}{}{}



\usepackage{authblk}

\title{A note on estimating Bass model parameters}
\author{Mengzhenyu Zhang, Hyun-Soo Ahn, Joline Uichanco\thanks{Stephen M. Ross School of Business, University of Michigan, Ann Arbor, MI, USA, \\ zhenyuzh, hsahn, jolineu@umich.edu}}

\date{August 2021}

\begin{document}

\maketitle

\section{Introduction}

\cite{bass1969new} proposed a model (the Bass model) for the timing of adoptions of a new product, where the adoption rate increases with the number of past adoptions. This model has been widely used in marketing and operations management literature to model new product demand over time. In this note, we provide a simple approach to estimate the Bass model parameters and prove the convergence rate.

\section{The model}

We define $(\Omega, \mathcal{F}, \mathbb{P}, \{\mathcal{F}\}_{t \ge 0})$ as a filtered probability space endowed with a cumulative adoption process $D = \{D_t, \ t \ge 0\}$ where $D_t$ is the cumulative adoptions by time $t$. Let $m$ be a positive integer that denotes the market size of potential customers. Hence, $D_t : \Omega \mapsto \{0, 1, \ldots, m\}$.  Since adoptions can only occur in unit increments, $D$ is a counting process.  
{Let $\{\mathcal{F}_t, \ t \geq 0\}$ be the history or filtration associated with the process of prices and adoptions, with $\mathcal F_t = \sigma((r_s, D_s), s \in [0,t])$.  We say that $\pi$ is a non-anticipating pricing policy if the price $r_t^\pi$ offered by $\pi$ at time $t$ is $\mathcal F_t$-measurable. } 
If customers are price-sensitive, a price change results in a change in the adoption rate.  To explicitly state the dependence in price, we will henceforth refer to the cumulative adoption as $D^\pi$ instead of $D$.  Without loss of generality, we assume that $D^\pi_0 = 0$ for any $\pi$, thus none of the consumers has purchased before time $t=0$. 

We denote the parameters of the Bass model as $\theta := (\alpha, \beta)$, {where $\alpha, \beta > 0$}. If at time $t$, the cumulative number of adoptions is $j$ and the seller sets price $r_t$, then the transition rate to the next $(j+1)$-st adoption is
\begin{align}
	\lambda(j, r_t) := \xi(j)\cdot x(r_t), \quad \mbox{for } j = 0, 1, \ldots, m, \label{eqn:lambda}
\end{align}
where
\begin{align}
    \xi(j) := (m - j) \left(\alpha + \beta \cdot \dfrac{j}{m} \right). \label{eqn:xi}
\end{align}
Note that $\xi(j)$ is the portion of the adoption rate unaffected by price.

\section{Estimation}

\cite{agrawal2021dynamic} proved that the expected mean squared estimation error of parameter $\beta$ cannot be better than $\Omega \left(m^2/n^3\right)$ given that $n$ is the number of adoptions (Lemma E.3). However, in what follows, we show that through a proper parameter transformation, the estimation error of $\beta$ under maximum likelihood estimation (MLE) is $\mathcal{O} \left(1/n\right)$.

We first conduct parameter transformation of the problem so that the variance of ML estimators does not grow infinitely as $d/m$ approaches zero and the second-order derivatives of the log-likelihood with respect to the unknown parameters do not interact with each other. We let $\alpha' := \alpha - \beta$ and $\beta' := \frac{\beta}{\alpha - \beta}$ so $\beta$ can be calculated from $\beta = \beta' \alpha'$ and $\alpha$ can be calculated from $\alpha = \alpha' + \beta' \alpha'$. Then, we have
\begin{equation*}
\begin{aligned}
    \ln \xi(d; \alpha, \beta) & = \ln (m-d) + \ln \left(\alpha + \beta \frac{d}{m}\right) \\
    & = \ln (m-d) + \ln\left(\alpha' \left(1 + (1 + \frac{d}{m})\beta'\right)\right) \\
    & = \ln (m-d) + \ln \alpha' + \ln\left(1 + (1 + \frac{d}{m})\beta'\right).
    \end{aligned}
\end{equation*}

We let both $\alpha', \beta'$ unknown and let $\hat \alpha'_t, \hat \beta'_t$ denote the estimated values from MLE at time $t$. The likelihood function is convenient to calculate under the Markovian Bass model.
We denote the continuously observed sequence of prices and cumulative sales at time $t$ as 
\begin{align}
	\widehat{\mathbf{U}}_t := \left\{ \left(\widehat{r}_s, \widehat{D}_s\right), \ 0 \leq s \leq t\right\}.
\end{align}
Since the adoption process follows a continuous-time Markov chain, the inter-adoption times are conditionally independent given the previous state information. Let $t_i$ be the time of the $i$th product adoption, where $i = 0, 1, 2, \ldots$ That is, at time $t_k$, the cumulative adoption is $ \widehat D_{t_k} = k$. The log-likelihood of $\widehat{\mathbf{U}}_t$ under a Markovian Bass model with parameters $\alpha', \beta'$ is
\begin{equation}
\begin{aligned}
\mathcal{L}_t (\widehat{\mathbf{U}}_t \mid \alpha',\beta') 
& = \sum_{i = 0}^{\widehat{D}_t - 1} \ln x(\widehat{r}_{t_{i + 1}}) + \sum_{i = 0}^{\widehat{D}_t - 1} \ln (m - i) + \sum_{i = 0}^{\widehat{D}_t - 1} \ln \alpha' + \sum_{i = 0}^{\widehat{D}_t - 1} \ln \left(1 + (1+ \frac{i}{m}) \beta' \right) \\
&\qquad  - \sum_{i = 0}^{\widehat{D}_t - 1}  \int_{t_{i}}^{t_{i + 1}} (m - i) \left(\alpha' \left(1 + (1 + \frac{i}{m})\beta'\right)\right) x(\widehat{r}_{s}) \mathrm{d} s\\
&\qquad  -  \int_{t_{\widehat{D}_t}}^{t} (m - \widehat{D}_t) \left(\alpha' \left(1 + (1 + \frac{\widehat{D}_t}{m})\beta'\right)\right) x(\widehat{r}_s) \mathrm{d} s. \\
\end{aligned}
\label{eqn:loglikelihood}
\end{equation}
Therefore, the ML estimators $\hat \alpha'_t, \hat \beta'_t$ are chosen such that the log-likelihood function $\mathcal{L}_t (\widehat{\mathbf{U}}_t \mid \alpha',\beta')$ is maximized at $\hat \alpha'_t, \hat \beta'_t$. The following result bounds the estimation error of $\hat \alpha'_t, \hat \beta'_t$.

\begin{theorem} \label{ParamErr}
We let $\alpha'_0,\beta'_0$ denote the true values. For any fixed time $t>0$ and $k \ge 2$, 
\begin{equation*} 
    \mathbb{E}_{\alpha'_0,\beta'_0} \left(\left(\hat{\alpha}'_{t} - \alpha'_0\right)^2 +  \left(\hat{\beta}'_{t} - \beta'_0\right)^2 \mid   D_t^\pi = n  \right) \leq \frac{\alpha_\theta}{n+1}, 
\end{equation*}
for some {$\alpha_\theta > 0$ that is independent of $m$, $t$ and $k$.} 
\end{theorem}

\begin{proof}
For simplicity of notation, we will use $D_t$ instead of $D_t^\pi$ to denote the cumulative adoptions at time $t$.

The ML estimators are finite since, from \eqref{eqn:loglikelihood},  if either $\hat \alpha'_t = +\infty$ or $\hat \beta'_t = +\infty$, then the likelihood function is 0. Hence, there exist finite $\bar \delta_1, \bar \delta_2$ such that $\hat \alpha'_t \le \alpha'_0 (1 + \bar \delta_1), \hat \beta'_t \le \beta'_0(1 + \bar \delta_2)$. 
Note that the ML estimator $\hat \theta_t = (\hat \alpha'_t, \hat \beta'_t)$ can be written as
\begin{equation*}
\hat{\theta}_t = \arg \max_{\theta \geq 0} \mathcal{L}_t (\widehat{\mathbf{U}}_t; \theta) = \theta_0 + \arg \min_{u \geq -\theta_0} - \sum_{i = 0}^{D_t} \ln \frac{f_i (\theta_0 + u)}{f_i (\theta_0)},
\end{equation*}
where $u = (u_{\alpha'}, u_{\beta'})$, $\theta_0 = (\alpha'_0,\beta'_0)$, and $f_i(\theta)$ is defined as follows:
{\small \begin{equation}\label{eqn:likelihood}
\begin{aligned}
	f_i(\theta) :=
	\begin{cases}
		(m - i) \left(\alpha' \left(1 + (1 + \frac{i}{m})\beta'\right)\right) x(\widehat{r}_{t_{i + 1}}) \exp\left(- (m - i) \left(\alpha' \left(1 + (1 + \frac{i}{m})\beta'\right)\right) \int_{t_{i}}^{t_{i + 1}} x(\widehat{r}_{s}) \mathrm{d} s \right), & \mbox{if } i = 0, 1,\ldots, D_t -1, \\
		\exp\left(- (m - D_t) \left(\alpha' \left(1 + (1 + \frac{D_t}{m})\beta'\right)\right) \int_{t_{D_t}}^{t} x(\widehat{r}_s) \mathrm{d} s\right), & \mbox{if } i= D_t.
	\end{cases}
\end{aligned}
\end{equation}}
If we denote the optimizer of the right-hand side as $\hat u = (\hat u_{\alpha'}, \hat u_{\beta'})$, then $\hat\theta_t = \theta_0 + \hat u$. 

We analyze the estimation error $|\hat \alpha'_t - \alpha'_0|$ first. 
Suppose $|\hat \alpha'_t - \alpha'_0| > \delta$ for some $\bar \delta_1 \alpha'_0 \ge \delta > 0$. This implies that $\hat u_{\alpha'}$ lies outside $[- \delta, \delta]$. Since the objective function on the right-hand-side is 0 when $u = 0$, and since the log-likelihood function is continuous and element-wise concave in $\alpha'$, then either
\begin{equation*}
- \sum_{i = 0}^{D_t} \ln \frac{f_i (\theta_0 + \delta e_1)}{f_i (\theta_0)} \leq 0 \hspace{5mm}\text{or}\hspace{5mm} -\sum_{i = 0}^{D_t} \ln \frac{f_i (\theta_0 - \delta e_1)}{f_i (\theta_0)} \leq 0,
\end{equation*}
where $e_1 := (1,0)$. Note that under the Markovian Bass model, the value $f_i(\theta)$ for any $\theta$ is stochastic since its value depends on $t_i$ and $t_{i+1}$, which are random adoption times. Here, $t_i$ denotes the time of the $i$-th adoption, where $i = 0, \ldots, D_t$.

Let $\mathbb{P}_{\theta_0}(\cdot)$ denote the  probability under a demand process that follows a Markovian Bass model with parameter vector $\theta_0 = (\alpha'_0, \beta'_0)$. 
Therefore,
\begin{align}
&\mathbb{P}_{\theta_0}  \left\{|\hat \alpha'_t - \alpha'_0| > \delta \right\} \notag \\
& \leq \mathbb{P}_{\theta_0} \left\{- \sum_{i = 0}^{D_t} \ln \frac{f_i (\theta_0+\delta e_1)}{f_i (\theta_0)} \leq 0 \right\} + \mathbb{P}_{\theta_0} \left\{- \sum_{i = 0}^{D_t} \ln \frac{f_i (\theta_0 - \delta e_1)}{f_i (\theta_0)} \leq 0 \right\} \notag \\
& \leq 2 \mathbb{P}_{\theta_0} \left\{- \sum_{i = 0}^{D_t} \ln \frac{f_i (\theta_0 + \delta e_1)}{f_i (\theta_0)} \leq 0  \right\}  = 2 \mathbb{P}_{\theta_0} \left\{\prod_{i = 0}^{D_t} \frac{f_i (\theta_0 + \delta e_1)}{f_i (\theta_0)} \geq 1  \right\} \notag\\
& = 2 \mathbb{P}_{\theta_0} \left\{\sqrt{\prod_{i = 0}^{D_t} \frac{f_i (\theta_0 + \delta e_1)}{f_i (\theta_0)}} \geq 1  \right\} \leq 2 \mathbb{E}_{\theta_0} \left(\sqrt{\prod_{i = 0}^{D_t} \frac{f_i (\theta_0 + \delta e_1)}{f_i (\theta_0)}} \right) \notag \\ 
& = 2  \mathbb{E}_{\theta_0} \left( \mathbb{E}_{\theta_0} \left(  \cdots \mathbb{E}_{\theta_0} \left( \mathbb{E}_{\theta_0} \left( \sqrt{\prod_{i = 0}^{D_t} \frac{f_i (\theta_0 + \delta e_1)}{f_i (\theta_0)}} \mid \mathcal F_{t_{D_t-1}} \right) \mid \mathcal F_{t_{D_t-2}} \right) \cdots \mid \mathcal F_{t_1} \right) \mid \mathcal F_0 \right) . \label{eqn:lemm3_prob_exp1}
\end{align}
The second inequality is because $f_i$ is an increasing function in $\alpha'$. The last equality is due to the law of iterated expectations. 

We next analyze \eqref{eqn:lemm3_prob_exp1} starting from the innermost conditional expectation. We have
\begin{equation} \label{eqn:tD-1}
\begin{aligned}
    & \mathbb{E}_{\theta_0} \left( \sqrt{\prod_{i = 0}^{D_t} \frac{f_i (\theta_0 + \delta e_1)}{f_i (\theta_0)}} \mid \mathcal F_{t_{D_t-1}} \right) = \sqrt{\prod_{i = 0}^{D_t-1} \frac{f_i (\theta_0 + \delta e_1)}{f_i (\theta_0)}} \mathbb{E}_{\theta_0} \left( \sqrt \frac{f_{D_t} (\theta_0 + \delta e_1)}{f_{D_t} (\theta_0)}\mid \mathcal F_{t_{D_t-1}}\right) \\
    & = \sqrt{\prod_{i = 0}^{D_t-1} \frac{f_i (\theta_0 + \delta e_1)}{f_i (\theta_0)}}\left( \int_{t_{D_t-1}}^{\infty}\sqrt \frac{f_{D_t} (\theta_0 + \delta e_1)}{f_{D_t} (\theta_0)} f_{D_t} (\theta_0) \mathrm d t_{D_t} \right)\\
    &= \sqrt{\prod_{i = 0}^{D_t-1} \frac{f_i (\theta_0 + \delta e_1)}{f_i (\theta_0)}} \left( \int_{t_{D_t-1}}^{\infty}\sqrt{f_{D_t} (\theta_0 + \delta e_1)}\sqrt{f_{D_t} (\theta_0)} \mathrm d t_{D_t}\right).
    \end{aligned}
\end{equation}
The first equality is because $\{f_i(\theta), i = 0, \ldots, D_t - 1\}$ are all $\mathcal F_{t_{D_t-1}}$-measurable. The second equality is because, given the information set $\mathcal F_{t_{D_t-1}}$, $f_{D_t}(\theta_0)$ is the conditional probability distribution of the adoption time $t_{D_t}$ under a Markovian Bass model with parameter $\theta_0$. Hence, we next want to derive a bound on $\int_{t_{D_t - 1}}^{\infty}\sqrt{f_{D_t} (\theta_0 + \delta e_1)}\sqrt{f_{D_t} (\theta_0)} \mathrm d t_{D_t}$.

Note that
\begin{align*}
    & \frac{1}{2} \int_{t_{D_{t}-1}}^\infty\left(\sqrt{f_{D_t} (\theta_0 + \delta e_1)} - \sqrt{f_{D_t} (\theta_0)} \right)^2 \mathrm d t_{D_t} \\
    & = \frac{1}{2}\int_{t_{D_t -1}}^\infty \left( f_{D_t}(\theta_0 + \delta e_1) + f_{D_t}(\theta_0) - 2 \sqrt{f_{D_t}(\theta_0 + \delta e_1) f_{D_t}(\theta_0)} \right) \mathrm d t_{D_t} \\
    & = 1 - \int_{t_{D_t - 1}}^\infty \sqrt{f_{D_t}(\theta_0 + \delta e_1) f_{D_t}(\theta_0)}\mathrm d t_{D_t},
\end{align*}
where the last equality is because the integral of the probability density function $\int_{t_{D_t-1}}^\infty f_{D_t}(\theta) \mathrm{d}t_{D_t}$ is equal to 1 for any $\theta$. Therefore,
\begin{equation} 
\begin{aligned}
    \int_{t_{D_t - 1}}^\infty \sqrt{f_{D_t}(\theta_0 + \delta e_1) f_{D_t}(\theta_0)}\mathrm d t_{D_t} & = 1 - \frac{1}{2} \int_{t_{D_{t}-1}}^\infty\left(\sqrt{f_{D_t} (\theta_0 + \delta e_1)} - \sqrt{f_{D_t} (\theta_0)} \right)^2 \mathrm d t_{D_t}.
    \end{aligned}
    \label{eqn:product_bound}
\end{equation}
The integral on the right-hand side is the Hellinger distance between $f_{D_t}(\theta_0 + \delta e_1)$ and $f_{D_t}(\theta_0)$, which are probability densities of the adoption time $t_{D_t}$.

Note that the Hellinger distance can be lower bounded by the K-L divergence (corollary 4.9 in \citealt{taneja2004relative}) provided the following condition holds. Specifically, 
\begin{equation} \label{eqn:H_KL}
\frac{1}{2} \int_{t_{D_{t}-1}}^\infty\left(\sqrt{f_{D_t} (\theta_0 + \delta e_1)} - \sqrt{f_{D_t} (\theta_0)} \right)^2 \mathrm d t_{D_t} \, \geq \, \frac{1}{4 \sqrt{R}} \mathbb{E}_{\theta_0} \left(\ln \frac{f_{D_t} (\theta_0)}{f_{D_t} (\theta_0 + \delta e_1)} \mid \mathcal{F}_{t_{D_t - 1}} \right),
\end{equation}
where $R$ is a constant such that $R \geq \max_{\delta \in [0, \bar \delta_1 \alpha'_0], t_{D_t}} \frac{1}{f_{D_t} (\theta_0 + \delta e_1)}$. Here, we can choose $R = 1/(\alpha'_0 + \alpha'_0\beta'_0)$ since  $\max_{\delta \in [0, \bar \delta_1 \alpha'_0], t_{D_t}} \frac{1}{f_{D_t} (\theta_0 + \delta e_1)} \leq 1/(m(\alpha'_0 + \alpha'_0\beta'_0)) \leq 1/(\alpha'_0 + \alpha'_0\beta'_0)$. Hence, with this choice, $R$ is independent of $m$ and of $t$. We will next bound the right-hand side of \eqref{eqn:H_KL}. 

Define $C_I := (\alpha'_0 (1 + \bar \delta_1))^2$. Note that
\begin{align*}
\frac{\partial^2  }{\partial \delta^2} \ln \frac{f_{D_t} (\theta_0)}{f_{D_t} (\theta_0 + \delta e_1)}  = \frac{1}{(\alpha'_0 + \delta)^2}\geq \frac{1}{(\alpha'_0 (1 + \bar \delta_1))^2} = \frac{1}{C_I},
\end{align*}
where the inequality is because $\alpha'_0 + \delta \le \alpha'_0 (1 + \bar \delta_1)$.

Furthermore, since the expectation of the Fisher score under the true parameter is zero, we have
\begin{align*} 
  \mathbb{E}_{\theta_0} \left(\frac{\partial}{\partial \delta} \ln \frac{f_{D_t} (\theta_0)}{f_{D_t} (\theta_0 + \delta e_1)} \Bigg|_{\delta = 0} \mid \mathcal{F}_{t_{D_t - 1}}  \right) = 0.
\end{align*}
Hence, we have
\begin{equation*}
\begin{aligned}
    & \mathbb{E}_{\theta_0} \left(\ln \frac{f_{D_t} (\theta_0)}{f_{D_t} (\theta_0 + \delta e_1)} \mid \mathcal{F}_{t_{D_t - 1}} \right) = \mathbb{E}_{\theta_0} \left(\int_0^{\delta} \frac{\partial}{\partial z} \ln \frac{f_{D_t} (\theta_0)}{f_{D_t} (\theta_0 + z e_1)} \mathrm d z\mid \mathcal{F}_{t_{D_t - 1}} \right) \\
    & = \mathbb{E}_{\theta_0} \left(\int_0^{\delta} \left(\frac{\partial}{\partial z} \ln \frac{f_{D_t} (\theta_0)}{f_{D_t} (\theta_0 + z e_1)} -  \frac{\partial}{\partial z} \ln \frac{f_{D_t} (\theta_0)}{f_{D_t} (\theta_0 + z e_1)} \Big|_{z = 0}\right)\mathrm d z\mid \mathcal{F}_{t_{D_t - 1}} \right) \\
    & = \mathbb{E}_{\theta_0} \left(\int_0^{\delta} \int_{0}^{z} \frac{\partial^2}{\partial {z'}^2} \ln \frac{f_{D_t} (\theta_0)}{f_{D_t} (\theta_0 + z' e_1)}\mathrm d z'\mid \mathcal{F}_{t_{D_t - 1}} \right) \ge \frac{1}{2C_I} \delta^2.
    \end{aligned}
\end{equation*}
Therefore, \eqref{eqn:H_KL} reduces to
\begin{align} \label{eqn:sqdiff_bd}
\frac{1}{4 \sqrt{R} C_I} \delta^2 & \leq \int_{t_{D_{t}-1}}^\infty \left(\sqrt{ f_{D_t} (\theta_0 + \delta e_1 )} - \sqrt{f_{D_t} (\theta_0)} \right)^2 \mathrm d t_{D_t}.
\end{align}

Hence, from \eqref{eqn:product_bound}, we have
\begin{equation*} 
\begin{aligned}
    &\int_{t_{D_t - 1}}^\infty \sqrt{f_{D_t}(\theta_0 + \delta e_1) f_{D_t}(\theta_0)}\mathrm d t_{D_t}  = 1 - \frac{1}{2} \int_{t_{D_{t}-1}}^\infty\left(\sqrt{f_{D_t} (\theta_0 + \delta e_1)} - \sqrt{f_{D_t} (\theta_0)} \right)^2 \mathrm d t_{D_t} \\
    &\hspace{10mm} \le \exp\left({- \frac{1}{2} \int_{t_{D_{t}-1}}^\infty\left(\sqrt{f_{D_t} (\theta_0 + \delta e_1)} - \sqrt{f_{D_t} (\theta_0)} \right)^2 \mathrm d t_{D_t} }\right) \leq \exp\left({- \frac{1}{8 \sqrt{R} C_I} \delta^2}\right),
    \end{aligned}
\end{equation*}
where the first inequality is because $e^{-x} \geq 1-x$ for all $x$. The second inequality is from \eqref{eqn:sqdiff_bd}. Hence, from \eqref{eqn:tD-1}, we have
\begin{align}
    \mathbb{E}_{\theta_0} \left( \sqrt{\prod_{i = 0}^{D_t} \frac{f_i (\theta_0 + \delta e_1)}{f_i (\theta_0)}} \mid \mathcal F_{t_{D_t-1}} \right) \leq 
    \sqrt{\prod_{i = 0}^{D_t-1} \frac{f_i (\theta_0 + \delta e_1)}{f_i (\theta_0)}} \cdot \exp \left({-\frac{1}{8\sqrt{R}C_I}\delta^2}\right). \label{eqn:innermost_bd}
\end{align}
This provides a bound for the innermost conditional expectation in \eqref{eqn:lemm3_prob_exp1}. 

Observe that all the terms in the bound \eqref{eqn:innermost_bd} are $\mathcal F_{t_{D_t-2}}$-measurable, except for the term 
\begin{equation*}
    \sqrt{f_{D_t-1}(\theta_0 + \delta e_1)/f_{D_t-1}(\theta_0)}.
\end{equation*}
Taking the conditional expectation of both sides in \eqref{eqn:innermost_bd} given $\mathcal F_{t_{D_t-2}}$, and using the same logic as the above arguments to bound the right-hand side, we have
\begin{align*}
    \mathbb{E}_{\theta_0} \left( \sqrt{\prod_{i = 0}^{D_t} \frac{f_i (\theta_0 + \delta e_1)}{f_i (\theta_0)}} \mid \mathcal F_{t_{D_t-2}} \right) \leq \sqrt{\prod_{i = 0}^{D_t-2} \frac{f_i (\theta_0 + \delta e_1)}{f_i (\theta_0)}} \cdot \exp \left( {-\frac{2}{8\sqrt{R}C_I}\delta^2} \right)
\end{align*}
We can proceed iteratively to evaluate \eqref{eqn:lemm3_prob_exp1} as we take conditional expectations given $\mathcal F_{t_{D_t}-3}$, $\mathcal F_{t_{D_t}-4}$, $\mathcal F_{0}$, resulting in
\begin{align*}
    \mathbb{E}_{\theta_0} \left( \sqrt{\prod_{i = 0}^{D_t} \frac{f_i (\theta_0 + \delta e_1)}{f_i (\theta_0)}} \right) \leq \mathbb{E}_{\theta_0} \left( \exp\left( - \frac{D_t + 1}{8 \sqrt{R} C_I} \delta^2 \right) \right)
\end{align*}

Hence, we have that
\begin{align*}
    \mathbb{P}_{\theta_0} \{|\hat{\alpha}'_t - \alpha'_0| > \delta \mid D_t = k\} \leq 2\mathbb{E}_{\theta_0} \left( \sqrt{\prod_{i = 0}^{D_t} \frac{f_i (\theta_0 + \delta e_1)}{f_i (\theta_0)}} \mid D_t = k \right) \leq 2\exp\left( - \frac{k + 1}{8 \sqrt{R} C_I} \delta^2 \right)
\end{align*}
if $\delta \le \bar\delta_1\alpha'_0$ and otherwise, $\mathbb{P}_{\theta_0} \{|\hat{\alpha}'_t - \alpha'_0| > \delta \mid D_t = k\} = 0$. 
This implies that
\begin{equation*}
\begin{aligned}
\mathbb{E}_{\theta_0} \left[(\hat \alpha'_t - \alpha'_0)^2 \mid D_t = n\right] &= \int_{0}^{\infty} \mathbb{P}_{\theta_0} \left\{(\hat \alpha'_t - \alpha'_0)^2 > \delta \mid D_t = n \right\} \, \mathrm{d} \delta  = \int_{0}^{\infty} \mathbb{P}_{\theta_0} \left\{|\hat \alpha'_t - \alpha'_0|^2 > \sqrt{\delta} \mid D_t = n \right\} \, \mathrm{d} \delta \\ 
& \leq \int_{0}^{\infty}  2\exp \left(- \frac{k + 1}{16 \sqrt{R} C_I} \delta \right) \mathrm{d} \delta = \frac{8 \sqrt{R} C_I}{n + 1} = \mathcal{O}\left( \frac{1}{n+1}\right).
\end{aligned}
\end{equation*}
Thus, we have that $\mathbb{E}_{\theta_0}\left[(\hat \alpha'_t - \alpha'_0)^2 \mid D_t = n\right] \leq \frac{\alpha_{\alpha'}}{n+1}$ 
where $\alpha_{\alpha'} := 8 \frac{\left(\alpha_0'(1+\bar \delta_1)\right)^2}{\sqrt{\alpha_0'+\alpha_0' \beta_0'}}$ is independent of $m$ and of $t$. 

Hence, to prove the lemma, we only need to show a similar bound for $\hat \beta'_t$. Similar bounds can be obtained for $\hat \beta'_t$ following the same steps with the only difference on the definition of $C_I$. 
We can safely write the second order derivative of the log-likelihood function with respect to $\beta'$. We have 
\begin{align*}
\mathbb{E}_{\theta_0} \left[\frac{\partial^2  }{\partial \delta^2} \ln \frac{f_{D_t} (\theta_0)}{f_{D_t} (\theta_0 + \delta e_2)} \mid \mathcal F_{t_{D_t-1}}  \right]  &= \mathbb{E}_{\theta_0} \left[ \frac{\left(1+\frac{D_t}{m}\right)^2}{\left(1 + \left(1+\frac{D_t}{m}\right)(\beta'_0 + \delta)\right)^2} \mid \mathcal F_{t_{D_t-1}} \right] \\
& \ge \mathbb{E}_{\theta_0} \left[ \frac{\left(1+\frac{D_t}{m}\right)^2}{\left(1 + \left(1+\frac{D_t}{m}\right)(\beta'_0(1 + \bar\delta_2))\right)^2} \mid \mathcal F_{t_{D_t-1}} \right],
\end{align*}
where the inequality is because $\beta'_0 + \delta \leq \beta'_0(1 + \bar \delta_2)$. Defining $C_I := \left(1 + \beta'_0 (1 + \bar \delta_2)\right)^2$, we have that
\begin{equation*}
    \mathbb{E}_{\theta_0} \left[ \frac{\left(1+\frac{D_t}{m}\right)^2}{\left(1 + \left(1+\frac{D_t}{m}\right)(\beta'_0(1 + \bar\delta_2))\right)^2} \mid \mathcal F_{t_{D_t-1}} \right] \ge \frac{1}{\left(1 + \beta_0' (1 + \bar \delta_2)\right)^2} = \frac{1}{C_{I}}.
\end{equation*}

Following the same steps in bounding the estimation error of $\hat \alpha'_t$, we know
\begin{align*}
    \mathbb{P}_{\theta_0} \{|\hat \beta'_t - \beta'_0| > \delta \mid D_t = n\} \leq 2\mathbb{E}_{\theta_0} \left( \sqrt{\prod_{i = 1}^{D_t} \frac{f_i (\theta_0 + \delta e_2)}{f_i (\theta_0)}} \mid D_t = n \right) \leq 2\exp\left( - \frac{k+1}{8 \sqrt{R}C_{I}}   \delta^2 \right).
\end{align*}

This implies that
\begin{equation*} \label{eqn:exp_error_q}
\begin{aligned}
\mathbb{E}_{\theta_0} \left[(\hat \beta'_t - \beta'_0)^2 \mid D_t = n\right] &= \int_{0}^{\infty} \mathbb{P}_{\theta_0} \left\{(\hat \beta'_t - \beta'_0)^2 > \delta \mid D_t = n \right\} \, \mathrm{d} \delta  = \int_{0}^{\infty} \mathbb{P}_{\theta_0} \left\{|\hat \beta'_t - \beta'_0|^2 > \sqrt{\delta} \mid D_t = n \right\} \, \mathrm{d} \delta \\ 
& \leq \int_{0}^{\infty}  2\exp \left(- \frac{n+1}{16 \sqrt{R}C_{I}}   \delta \right) \mathrm{d} \delta = \frac{8 \sqrt{R} C_I}{n+1} = \frac{8 \sqrt{R} \left(1 + \beta'_0 (1 + \bar \delta_2)\right)^2}{n+1}. 
\end{aligned}
\end{equation*}
Thus, we have that $\mathbb{E}_{\theta_0}\left[ {(\hat \beta'_t - \beta'_0)^2} \mid D_t = n\right] \leq \frac{\alpha_{\beta'}}{n+1}$ where $\alpha_{\beta'} := 8{\sqrt{R}}\left(1 + \beta'_0 (1 + \bar \delta_2)\right)^2$ is independent of $m$ and of $t$.

Hence, we prove the lemma with $\alpha_\theta := 8 {\sqrt{R}}\max \left\{ \left((1 + \bar \delta_1)\alpha'_0\right)^2, \left(1 + \beta'_0 (1 + \bar \delta_2)\right)^2\right\}$, and $R = 1/(\alpha'_0+\alpha'_0+\beta'_0)$. Note that $\alpha_\theta$ does not depend on $t$ and $m$.

\end{proof}

Note that we can also analyze the variance of $\hat \beta'_t$ under MLE. To do that, we need to analyze the bounds of the inverse of the fisher information:
\begin{equation*}
     \frac{1}{ -\sum_{d =1}^n\frac{\partial^2  }{\partial \beta'^2} \ln \xi(d; \alpha', \beta')} = 1 / \sum_{d=1}^n\left(\frac{(1 + d/m)^2}{\left(1 + (1 + d/m)\beta'\right)^2}\right).
\end{equation*}
We know that 
\begin{equation*}
    \sum_{d=1}^n \frac{1^2}{\left(1 + (1 + 1)\beta'\right)^2}\le \sum_{d=1}^n\left(\frac{(1 + d/m)^2}{\left(1 + (1 + d/m)\beta'\right)^2}\right) \le \sum_{d=1}^n \frac{2^2}{\left(1 + \beta'\right)^2}.
\end{equation*}
Thus, the variance of $\hat \beta'_t$ is upper bounded by 
\begin{equation*}
    \frac{1}{ -\sum_{d =1}^n\frac{\partial^2  }{\partial \beta'^2} \ln \xi(d; \alpha', \beta')} \le 1/\sum_{d=1}^n \frac{1^2}{\left(1 + (1 + 1)\beta'\right)^2} = \mathcal{O}(1/n).
\end{equation*}
This is consistent with the bound of mean squared error derived in \cref{ParamErr}.

We note that transforming parameters from $\alpha,\beta$ to $\alpha',\beta'$ is the same as pre-processing the data. Pre-processing the data could help with the estimation. For example, by our transformation of the parameters, we are actually viewing $1 + D_t/m$ as the data multiplied to $\beta'$, which avoids the issue of the large estimation variance of $\hat \beta_t$ if we directly estimate $\beta$ from MLE. This is a common practice in parameter estimation when data values tend to be small. 

\bibliographystyle{plainnat}
\bibliography{ref}
\end{document}